# MATURITY RANDOMIZATION FOR STOCHASTIC CONTROL PROBLEMS

By Bruno Bouchard, Nicole El Karoui and Nizar Touzi

*Université Paris VI, École Polytechnique and Université Paris IX*

We study a maturity randomization technique for approximating optimal control problems. The algorithm is based on a sequence of control problems with random terminal horizon which converges to the original one. This is a generalization of the so-called *Canadization* procedure suggested by Carr [*Review of Financial Studies* **II** (1998) 597–626] for the fast computation of American put option prices. In addition to the original application of this technique to optimal stopping problems, we provide an application to another problem in finance, namely the super-replication problem under stochastic volatility, and we show that the approximating value functions can be computed explicitly.

**1. Introduction.** It is well known that the arbitrage-free price of an American put in a complete market is the value of an optimal stopping problem, which corresponds in a Markov framework to a free boundary problem. For a finite horizon, no explicit formula for this value is known in general. An explicit solution does exist in the infinite-horizon case when the reward process is defined by a Lévy process; see, for example, [10].

The maturity randomization technique introduced by Carr [2] provides an interesting algorithm for the computation of a finite-horizon optimal stopping problem by passing to a sequence of infinite-horizon stopping problems. This technique is well established in the literature, and is referred to as the *Canadization* procedure; see, for example, [9]. We shall review this technique in Section 2.

However, the original paper of Carr [2] does not report a proof of consistency of this technique. Instead, there is an intuitive discussion of the theoretical foundations of the algorithm through appeal to the dynamic programming principle. Although this argument seems to be very intuitive, it does not apply to this particular context, as the random times introduced in the maturity randomization algorithm are independent of the filtration relative to the class of stopping times. The numerical evidence provided in [2] shows the excellent performance of this method.

In this paper we extend this approach to general finite-horizon stochastic control problems, including optimal stopping problems. The consistency of the algorithm is proved in this general framework. These results are contained in Section 3, and the application to optimal stopping problems is reported in Section 4. We conclude the paper by studying an example of stochastic control problem from finance, namely the problem of hedging in the uncertain volatility model. The value function of this problem can be characterized as the unique solution of a fully nonlinear partial differential equation. Applying the maturity randomization technique in this context, we are reduced to a sequence of nonlinear ordinary differential equations that can be solved explicitly.

**2. Solving the American put problem by maturity randomization.** In this section we review the numerical procedure suggested by Carr [2] for a fast numerical computation of the American put price. Let $(\Omega, \mathcal{F}, \mathbb{P})$ be a complete probability space supporting a real-valued Brownian motion $W = \{W(t), t \geq 0\}$. We denote by $\mathbb{F} = \{\mathcal{F}_t,\ t \geq 0\}$ the $\mathbb{P}$-completion of the canonical filtration of $W$.

For every $t \geq 0$, the set $\mathcal{T}_t(\mathbb{F})$ is the collection of all $\mathbb{F}$-stopping times $\tau \geq t$ $\mathbb{P}$-a.s.

2.1. *The American put problem.* Let $S$ be the process defined by

$$S(t) = S(0) \exp\left[\left(r - \frac{\sigma^2}{2}\right)t + \sigma W(t)\right], \qquad t \geq 0,$$

where $S(0)$ is some given initial data, and $r, \sigma > 0$ are given parameters. The main purpose of [2] is to compute the value of the following optimal stopping problem:

$$(2.1) \qquad V_0 := \sup_{\tau \in \mathcal{T}_0(\mathbb{F})} \mathbb{E}[e^{-r(\tau \wedge T)} g(S(\tau \wedge T))],$$

where $T > 0$ is some given finite horizon, and

$$g(x) := [K - x]^+ \qquad \text{for some positive constant } K.$$

We introduce the so-called Snell envelope of the reward process $\{e^{-r(t \wedge T)} g(S(t \wedge T)),\ t \geq 0\}$:

$$V(t) := \operatorname*{ess\,sup}_{\tau \in \mathcal{T}_t(\mathbb{F})} \mathbb{E}[e^{-r(\tau \wedge T)} g(S(\tau \wedge T)) | \mathcal{F}_t],$$

whose analysis provides a complete characterization of the solution of (2.1). From the Markov property of the process $S$, the above Snell envelope can be written as

$$V(t) = v(t, S(t)) \qquad \text{for all } t \geq 0,$$

where $v$ is the value function of the dynamic version of the optimal stopping problem (2.1):

$$v(t,x) := \sup_{\tau \in \mathcal{T}_t(\mathbb{F})} \mathbb{E}[e^{-r(\tau \wedge T)} g(S(\tau \wedge T)) | S(t) = x].$$

2.2. *Maturity randomization.* The main idea of [2] is to reduce the problem of computation of $V_0$ to a sequence of infinite-horizon optimal stopping problems, which are well known to be easier to solve. Indeed when $T = +\infty$, it follows from the homogeneity of the process $S$ that the dependence of the value function $v$ on the time variable is given by

$$v(t,x) = e^{rt} v(0,x) \qquad \text{for all } (t,x) \in \mathbb{R}_+^2,$$

and the problem reduces to finding the dependence of $v$ on the $x$ variable. In many instances, this dependence can be found explicitly. We now describe Carr's procedure in different steps.

*Step* 1. A sequence of infinite-horizon optimal stopping problems is created by approximating the fixed finite maturity $T$ by a sequence of random variables. Let $(\xi^k)_{k \geq 0}$ be a sequence of random variables satisfying the following requirements:

(2.2) $\qquad \xi^k$ are i.i.d. nonnegative random variables with $\mathbb{E}[\xi^k] = 1$,

(2.3) $\qquad \xi^k$ is independent of $\mathbb{F}$ for every $k \geq 0$.

By the law of large numbers, it follows from (2.2) that

$$T_n^n := \frac{T}{n} \sum_{j=1}^{n} \xi^j \longrightarrow T, \qquad \mathbb{P}\text{-a.s.}$$

It is then natural to introduce the approximation

$$v_n(x) := \sup_{\tau \in \mathcal{T}_0(\mathbb{F})} \mathbb{E}[e^{-r(\tau \wedge T_n^n)} g(S(\tau \wedge T_n^n)) | S(0) = x].$$

In the sequel, we shall need the extended notation

$$v_n^k(x) := \sup_{\tau \in \mathcal{T}_0(\mathbb{F})} \mathbb{E}[e^{-r(\tau \wedge T_n^k)} g(S(\tau \wedge T_n^k)) | S(0) = x],$$

where

$$T_n^k := \frac{T}{n} \sum_{j=1}^{k} \xi^j \qquad \text{for } k \leq n,$$

and we observe that $v_n^n = v_n$.

*Step* 2. We next observe that
$$T_n^k = T_n^{k-1} + \zeta_n^k \qquad \text{where } \zeta_n^k := \frac{T}{n}\xi^k,$$
and we use property (2.3) of the random variables $(\xi^j)$ to write
$$\begin{aligned}v_n^k(x) &= \sup_{\tau\in\mathcal{T}_0(\mathbb{F})} \mathbb{E}[e^{-r(\tau\wedge T_n^k)}g(S(\tau\wedge T_n^k))\mathbb{1}_{\{\tau\leq\zeta_n^k\}}\\ &\qquad\qquad + e^{-r(\tau\wedge T_n^k)}g(S(\tau\wedge T_n^k))\mathbb{1}_{\{\tau>\zeta_n^k\}}|S(0)=x]\\ &= \sup_{\tau\in\mathcal{T}_0(\mathbb{F})} \mathbb{E}[e^{-r\tau}g(S(\tau))G_n^k(\tau)\\ &\qquad\qquad + e^{-r(\tau\wedge T_n^k)}g(S(\tau\wedge T_n^k))\mathbb{1}_{\{\tau>\zeta_n^k\}}|S(0)=x],\end{aligned}$$
where
$$G_n^k(t) := \mathbb{P}[\zeta_n^k \geq t].$$

*Step* 3. By a formal argument, Carr claims that the latter supremum can be written as
$$\begin{aligned}(2.4)\qquad v_n^k(x) &= \sup_{\tau\in\mathcal{T}_0(\mathbb{F})} \mathbb{E}[e^{-r\tau}g(S(\tau))G_n^k(\tau)\\ &\qquad\qquad + e^{-r\zeta_n^k}v_n^{k-1}(S(\zeta_n^k))\mathbb{1}_{\{\tau>\zeta_n^k\}}|S(0)=x].\end{aligned}$$
Let us point out that Carr fully recognizes that he is not providing a rigorous proof for the convergence of the scheme. We shall elaborate further on this point later on, but let us only observe that, at a first glance, this equality seems to follows from

  (i) the classical *dynamic programming principle*,
  (ii) the homogeneous feature of the problem.

*Step* 4. Using again the fact that $\zeta_n^k$ is independent of $\mathbb{F}$, the above formula (2.4) can be written as
$$v_n^k(x) = \sup_{\tau\in\mathcal{T}_0(\mathbb{F})} \mathbb{E}\left[e^{-r\tau}g(S(\tau))G_n^k(\tau) - \int_0^\tau e^{-rt}v_n^{k-1}(S(t))\,dG_n^k(t)\Big|S(0)=x\right].$$
Finally, fix the distribution of $\xi^i$ to be exponential with unit parameter. Then
$$G_n^k(t) = e^{-nt/T} \qquad \text{for all } t\geq 0,$$
and we obtain the following recursive formula:
$$(2.5)\quad v_n^k(x) = \sup_{\tau\in\mathcal{T}_0(\mathbb{F})} \mathbb{E}\left[e^{-r_n\tau}g(S(\tau)) + \frac{n}{T}\int_0^\tau e^{-r_n t}v_n^{k-1}(S(t))\,dt\Big|S(0)=x\right],$$

where we defined the parameter

$$r_n := r + \frac{n}{T}.$$

*Step* 5. In the case of the American put option, Carr was able to write a beautiful explicit formula which relates $v_n^k$ to $v_n^{k-1}$; that is, given the function $v_n^{k-1}$, the optimal stopping problem (2.5) is solved explicitly. Together with the use of the Richardson extrapolation technique, this produces a fast and accurate approximation of the American put option value.

2.3. *Consistency and extension to general control problems.* The first objective of this paper was to provide a rigorous proof of consistency for the scheme described in the previous paragraph. This opened the door for a much larger generality of this technique.

Our first attempt for the proof of consistency is to justify the crucial equality (2.4). Unfortunately, the dynamic programming principle does not apply in this context, as $\zeta_n^k$ is independent of the filtration $\mathbb{F}$. Our first main result is that, although this equality *may not hold*, the scheme suggested by Carr by the recursive formula (2.5) is consistent. The proof is provided in Section 4.2.

In Section 4 the above result is established for general optimal stopping problems, thus dropping the Markov and the homogeneity assumptions on the reward process. The random variables $\xi^k$ are also allowed to have different distributions. This could be exploited as an error reduction factor. We leave this point for further research.

In Section 3 we prove that the maturity randomization technique applies to general stochastic control problems, and mixed stopping/control problems.

We conclude the paper by providing another interesting example where the maturity randomization technique leads to an explicit recursive relation. The example studied in Section 5 consists in the problem of hedging a European contingent claim in the context of the uncertain volatility model, that is, the diffusion coefficient is only known to lie in between two bounds.

## 3. Approximating control problems by maturity randomization.

3.1. *The control problems.* We now consider a general probability space $(\Omega, \mathcal{A}, \mathbb{P})$ endowed with a filtration $\mathbb{F} = \{\mathcal{F}_t\}_{t \geq 0}$ satisfying the usual conditions, and we assume that $\mathcal{F}_0$ is trivial. Importantly, we do not assume that $\mathcal{A} = \mathcal{F}_\infty$ in order to allow for other sources of randomness.

Given a set $U$ of (deterministic) functions from $\mathbb{R}_+$ to $\mathbb{R}^d$, $d \geq 1$, we denote by $\tilde{\mathcal{U}}(\mathbb{F})$ the collection of all $\mathbb{F}$-adapted processes $\nu$ such that

$$t \longmapsto \nu(t, \omega) \in U \qquad \text{for almost every } \omega \in \Omega.$$

The controlled state process is defined by a map

$$\nu \in \tilde{\mathcal{U}}(\mathbb{F}) \longmapsto Y^\nu \in L^0_{\mathbb{F}}(\mathbb{R}_+ \times \Omega, \mathbb{R}),$$

where $L^0_{\mathbb{F}}(\mathbb{R}_+ \times \Omega, \mathbb{R})$ is the set of all $\mathbb{F}$-progressively measurable processes valued in $\mathbb{R}$, and

(3.1) $\qquad\qquad Y^\nu(0) =: Y(0) \qquad$ is independent of $\nu$.

The set $\mathcal{U}(\mathbb{F})$ of $\mathbb{F}$-admissible control processes is a subset of the collection of elements $\nu \in \tilde{\mathcal{U}}(\mathbb{F})$. We assume that this set of controls is stable under bifurcation at deterministic times, that is,

(H$\mathcal{U}$)  For all $\nu_1, \nu_2 \in \mathcal{U}(\mathbb{F}), t \geq 0$ and $A \in \mathcal{F}_t$,

$$\nu_1 = \nu_2 \text{ on } [0, t) \text{ $\mathbb{P}$-a.s.} \quad \Longrightarrow \quad \nu_1|t_A|\nu_2 := \nu_1 \mathbb{1}_A + \nu_2 \mathbb{1}_{A^c} \in \mathcal{U}(\mathbb{F}).$$

Notice that this condition is slightly weaker than the stability by bifurcation at stopping times introduced in [5].

REMARK 3.1.  Assumption (H$\mathcal{U}$) is weaker than the usual stability under concatenation property:

(H$\mathcal{U}$)$'$ For all $\nu_1, \nu_2 \in \mathcal{U}(\mathbb{F})$ and $\tau \in \mathcal{T}_0(\mathbb{F})$, $\nu_1 \mathbb{1}_{[0,\tau)} + \nu_2 \mathbb{1}_{[\tau,\infty)} \in \mathcal{U}(\mathbb{F})$,

which is not satisfied for the optimal stopping problems studied in Section 4. In Section 3.3, we shall use a weak version of (H$\mathcal{U}$)$'$:

(H$\mathcal{U}$)$''$ For all $\nu_1, \nu_2 \in \mathcal{U}(\mathbb{F})$ and $t \geq 0$, $\nu_1 \mathbb{1}_{[0,t)} + \nu_2 \mathbb{1}_{[t,\infty)} \in \mathcal{U}(\mathbb{F})$.

We are interested in computing

(3.2) $$\sup_{\nu \in \mathcal{U}(\mathbb{F})} \mathbb{E}[Y^\nu(T)].$$

Following the maturity randomization technique of [2], we introduce a sequence of approximating control problems. We denote by $\mathcal{I}^{\mathbb{F}}$ the collection of all nonnegative random variables $\xi$ which are independent of $\mathcal{F}_\infty$, that is,

$$\mathbb{E}[\xi \mathbb{1}_A] = \mathbb{P}[A]\mathbb{E}[\xi] \qquad \text{for any } A \in \mathcal{F}_\infty.$$

Given some integer $n \geq 1$, we next consider a sequence $(\zeta^j)_{1 \leq j \leq n}$ of independent random variables in $\mathcal{I}^{\mathbb{F}}$, and we set

$$T^k := \sum_{j=1}^k \zeta^j \qquad \text{for } 1 \leq k \leq n, \qquad T^0 := 0.$$

We denote by $m$ the law of $(\zeta^1, \ldots, \zeta^n)$ under $\mathbb{P}$, that is,

$$m(A_1 \times \cdots \times A_n) = \mathbb{P}[\zeta^1 \in A_1, \ldots, \zeta^n \in A_n]$$
$$= \prod_{j=1}^n m^j(A_j) \qquad \text{for all } A_1, \ldots, A_n \in \mathcal{B}_{\mathbb{R}_+},$$

where $\mathcal{B}_{\mathbb{R}_+}$ denotes the Borel tribe of $\mathbb{R}_+$, and $m^j$ denotes the law of $\zeta^j$.

The maturity randomization algorithm is defined as follows:

$$V_0^\nu = Y^\nu, \qquad \nu \in \mathcal{U}(\mathbb{F}), \tag{3.3}$$

and for $k = 0, \ldots, n-1$

$$V_{k+1}^\nu(t) = \operatorname*{ess-sup}_{\mu \in \mathcal{U}(\mathbb{F};t,\nu)} \mathbb{E}[\bar{V}_k^\mu(t + \zeta^{n-k})|\mathcal{F}_t], \qquad t \geq 0, \tag{3.4}$$

where $\bar{V}_k^\mu$ is an $(\Omega \times \mathbb{R}_+, \mathcal{F} \otimes \mathcal{B}_{\mathbb{R}_+})$-measurable aggregating process for $V_k^\mu$ [see assumption (HV) below], and

$$\mathcal{U}(\mathbb{F};t,\nu) := \{\mu \in \mathcal{U}(\mathbb{F}) : \mu = \nu \text{ on } [0,t) \text{ } \mathbb{P}\text{-a.s.}\}.$$

In order to give a sense to the above expressions, we assume that

(HY) There is a uniformly integrable martingale $M^Y$ such that, for each $\nu \in \mathcal{U}(\mathbb{F})$, $|Y^\nu(t)| \leq M^Y(t)$ for all $t \geq 0$ $\mathbb{P}$-a.s.
(HV) For each $\nu \in \mathcal{U}(\mathbb{F})$ and $1 \leq k \leq n-1$, there is an $(\Omega \times \mathbb{R}_+, \mathcal{F} \otimes \mathcal{B}_{\mathbb{R}_+})$-measurable process $\bar{V}_k^\nu$ such that $\bar{V}_k^\nu(t) = V_k^\nu(t)$ $\mathbb{P}$-a.s. for all $t \geq 0$.

REMARK 3.2. Assumption (HY) implies that, for each $\nu \in \mathcal{U}(\mathbb{F})$ and $0 \leq k \leq n$, $|V_k^\nu(t)| \leq M^Y(t)$ $\mathbb{P}$-a.s. for all $t \geq 0$. Indeed, assume that the assertion is true for some $0 \leq k < n$. Since $\zeta^{n-k}$ is independent of $\mathcal{F}$, using Fubini's lemma in (3.4) leads to

$$\begin{aligned}
\bar{V}_{k+1}^\nu(t) &\leq \operatorname*{ess-sup}_{\mu \in \mathcal{U}(\mathbb{F};t,\nu)} \mathbb{E}[|\bar{V}_k^\mu(t + \zeta^{n-k})||\mathcal{F}_t] \\
&= \operatorname*{ess-sup}_{\mu \in \mathcal{U}(\mathbb{F};t,\nu)} \int \mathbb{E}[|\bar{V}_k^\mu(t + z^{n-k})||\mathcal{F}_t] m(dz) \\
&\leq \int \mathbb{E}[M^Y(t + z^{n-k})|\mathcal{F}_t] m(dz) \\
&= M^Y(t), \qquad \mathbb{P}\text{-a.s.}
\end{aligned}$$

The same argument also shows that the expectations in (3.4) are well defined.

REMARK 3.3. (i) Assumption (HV) is necessary since $V_k^\mu(t + \zeta^{n-k})$ may not be defined as a measurable map from $\Omega$ into $\mathbb{R}$.

(ii) Observe that $\bar{V}_0^\nu = V_0^\nu$ from the conditions on the controlled process $Y^\nu$.

(iii) In the usual literature on stochastic control (see, e.g., [5]), (3.4) is shown to define a supermartingale family, as a consequence of the stability by bifurcation property. This is the key point in order to prove the existence of a ladlag aggregating supermartingale. Unfortunately, these results do not apply in our framework. Indeed, the time $t$ appears on the right-hand side of (3.4) both in the controlled process and in the conditioning, so that the problem (3.4) does not fit in the classical class of stochastic control problems.

(iv) In Sections 3.3 and 4.2 we shall provide sufficient conditions ensuring the existence of a ladlag modification of $V_k^\nu$. This will be obtained by assuming an exponential distribution for the maturity randomizing random variables $\zeta^k$, and observing that (3.4) reduces, in this case, to a classical stochastic control problem.

REMARK 3.4. For later use, notice that, under assumption (H$\mathcal{U}$), for $t_2 \geq t_1 \geq 0$ and $\nu \in \mathcal{U}(\mathbb{F})$

$$\mathcal{U}(\mathbb{F}; t_1, \nu) \supset \{\nu_2 \in \mathcal{U}(\mathbb{F}; t_2, \nu_1), \nu_1 \in \mathcal{U}(\mathbb{F}; t_1, \nu)\}.$$

Since $\mathcal{U}(\mathbb{F}; 0, \nu) = \mathcal{U}(\mathbb{F})$, we shall simply write

(3.5) $$V_k(0) := V_k^\nu(0) \qquad \text{for } k \leq n.$$

3.2. *The convergence result.* We start with the following easy lemma which will be used later to derive an upper bound for $V_n(0)$.

LEMMA 3.1. *Under assumptions* (HY), (HV) *and* (H$\mathcal{U}$), *for all $k \geq 0$, $t \geq 0$, $\nu \in \mathcal{U}(\mathbb{F})$*,

$$\mathbb{E}\left[\operatorname*{ess\,sup}_{\mu \in \mathcal{U}(\mathbb{F}; t; \nu)} \mathbb{E}[\bar{V}_k^\mu(t + \zeta^{n-k}) | \mathcal{F}_t]\right] = \sup_{\mu \in \mathcal{U}(\mathbb{F}; t; \nu)} \mathbb{E}[\bar{V}_k^\mu(t + \zeta^{n-k})].$$

*In particular,*

$$\mathbb{E}[\bar{V}_{k+1}^\nu(t)] = \sup_{\mu \in \mathcal{U}(\mathbb{F}; t; \nu)} \mathbb{E}[\bar{V}_k^\mu(t + \zeta^{n-k})].$$

PROOF. Under assumption (H$\mathcal{U}$), the family

$$\{\mathbb{E}[\bar{V}_k^\mu(t + \zeta^{n-k}) | \mathcal{F}_t], \mu \in \mathcal{U}(\mathbb{F}; t; \nu)\}$$

is directed upward. We can then find a sequence $\mu_j \in \mathcal{U}(\mathbb{F}; t; \nu)$ such that

$$\operatorname*{ess\,sup}_{\mu \in \mathcal{U}(\mathbb{F}; t; \nu)} \mathbb{E}[\bar{V}_k^\mu(t + \zeta^{n-k}) | \mathcal{F}_t] = \lim_{j \to \infty} \uparrow \mathbb{E}[\bar{V}_k^{\mu_j}(t + \zeta^{n-k}) | \mathcal{F}_t], \qquad \mathbb{P}\text{-a.s.}$$

By the monotone convergence theorem, this implies that

$$\mathbb{E}\left[\operatorname*{ess\,sup}_{\mu \in \mathcal{U}(\mathbb{F}; t; \nu)} \mathbb{E}[\bar{V}_k^\mu(t + \zeta^{n-k}) | \mathcal{F}_t]\right] \leq \sup_{\mu \in \mathcal{U}(\mathbb{F}; t; \nu)} \mathbb{E}[\bar{V}_k^\mu(t + \zeta^{n-k})].$$

The converse inequality is obviously satisfied. The second statement of the lemma then follows from the definition of $V_{k+1}^\nu$ in (3.4). □

We are now ready for the main result of this section.

THEOREM 3.1. *Let assumptions* (HY), (HV) *and* (H$\mathcal{U}$) *hold. Then*

$$(3.6) \quad \sup_{\nu \in \mathcal{U}(\mathbb{F})} \mathbb{E}[Y^\nu(T^n)] \leq V_n(0) \leq \int \sup_{\nu \in \mathcal{U}(\mathbb{F})} \mathbb{E}[Y^\nu(z^1 + \cdots + z^n)] m(dz).$$

PROOF. 1. We first prove the upper bound. Fix $1 \leq k \leq n-1$, $\nu \in \mathcal{U}(\mathbb{F})$ and $t \geq 0$. Since $\zeta^{n-k}$ is independent of $\mathcal{F}_\infty$, it follows from assumption (HY) and Remark 3.2 that we can use Fubini's lemma to get

$$\mathbb{E}[\bar{V}_k^\nu(t + \zeta^{n-k})] = \int \mathbb{E}[\bar{V}_k^\nu(t + z^{n-k})] m(dz),$$

where we use the notation $z = (z^1, \ldots, z^n)$. By Lemma 3.1, this can be written as

$$\mathbb{E}[\bar{V}_k^\nu(t + \zeta^{n-k})] = \int \sup_{\mu \in \mathcal{U}(\mathbb{F}; t+z^{n-k}; \nu)} \mathbb{E}[\bar{V}_{k-1}^\mu(t + z^{n-k} + \zeta^{n-k+1})] m(dz).$$

In view of Remark 3.4, the upper bound of Theorem 3.1 then follows from an easy induction.

2. In order to provide the lower bound, we first show that for all $\nu \in \mathcal{U}(\mathbb{F})$:

$$(3.7) \quad \begin{aligned} \mathbb{E}[\bar{V}_k^\nu(T^{n-k})] &= \mathbb{E}[\bar{V}_k^\nu(\zeta^{n-k} + T^{n-k-1})] \\ &\leq \mathbb{E}[\bar{V}_{k+1}^\nu(T^{n-k-1})], \qquad k \leq n-1. \end{aligned}$$

Indeed, since $(\zeta^k)_{k \leq n}$ are independent random variables in $\mathcal{I}^\mathbb{F}$, we have

$$\mathbb{E}[\bar{V}_k^\nu(\zeta^{n-k} + T^{n-k-1})] = \mathbb{E}\left[\int_0^\infty \bar{V}_k^\nu(\zeta^{n-k} + t) \, dF(t)\right],$$

where $F(t) := \mathbb{P}[T^{n-k-1} \leq t]$ is the cumulative probability distribution of $T^{n-k-1}$. We next use Fubini's lemma together with the definition of $V_k^\nu$ in (3.4) to obtain

$$\begin{aligned} \mathbb{E}[\bar{V}_k^\nu(\zeta^{n-k} + T^{n-k-1})] &= \int_0^\infty \mathbb{E}[\mathbb{E}[\bar{V}_k^\nu(\zeta^{n-k} + t) | \mathcal{F}_t]] \, dF(t) \\ &\leq \int_0^\infty \mathbb{E}[V_{k+1}^\nu(t)] \, dF(t) \\ &= \int_0^\infty \mathbb{E}[\bar{V}_{k+1}^\nu(t)] \, dF(t) \\ &= \mathbb{E}[\bar{V}_{k+1}^\nu(T^{n-k-1})]. \end{aligned}$$

By (3.3), (3.5), it follows by using repeatedly (3.7) that

$$\mathbb{E}[Y^\nu(T^n)] = \mathbb{E}[\bar{V}_0^\nu(T^n)] \leq \mathbb{E}[V_n^\nu(0)] = V_n(0).$$

Since $\nu$ is an arbitrary control in $\mathcal{U}(\mathbb{F})$, this provides the lower bound announced in Theorem 3.1. □

We now consider sequences $\{(\zeta_n^k)_{k \leq n}\}_{n \geq 1}$ of random variables in $\mathcal{I}^{\mathbb{F}}$. We define the corresponding sequence $\{(V_k^{\bar{\nu},n})_{k \leq n}\}_{n \geq 1}$, where, for each $n$, $(V_k^{\nu,n})_{k \leq n}$ is defined as in (3.3), (3.4) with the sequence $(\zeta_n^k)_{k \leq n}$. For each $n \geq 1$, we define

$$T_n^n := \sum_{j=1}^{n} \zeta_n^j,$$

and we denote by $m_n$ the law of $(\zeta_n^1, \ldots, \zeta_n^n)$. Using the bounds of Theorem 3.1, we can provide conditions under which $V_n^n(0) = V_n^{\nu,n}(0)$ converges to the value of the initial control problem (3.2).

COROLLARY 3.1. *Let assumptions* (HY), (HV) *and* (H$\mathcal{U}$) *hold, and suppose that the triangular array* $(\zeta_n^k)$ *satisfies*

$$T_n^n \longrightarrow T \in (0, \infty) \qquad \text{in probability.}$$

(i) *Assume that*

(3.8) $\quad t > 0 \longmapsto \mathbb{E}[Y^\nu(t)] \qquad$ *is continuous at $t = T$ for all $\nu \in \mathcal{U}(\mathbb{F})$.*

*Then* $\liminf_{n \to \infty} \mathbb{E}[Y^\nu(T_n^n)] \geq \mathbb{E}[Y^\nu(T)]$ *for all $\nu \in \mathcal{U}(\mathbb{F})$.*

(ii) *Assume that*

(3.9) $\quad t > 0 \longmapsto \sup_{\nu \in \mathcal{U}(\mathbb{F})} \mathbb{E}[Y^\nu(t)] \qquad$ *is continuous at $t = T$.*

*Then* $\limsup_{n \to \infty} \int \sup_{\nu \in \mathcal{U}(\mathbb{F})} \mathbb{E}[Y^\nu(z^1 + \cdots + z^n)] m_n(dz) \leq \sup_{\nu \in \mathcal{U}(\mathbb{F})} \mathbb{E}[Y^\nu(T)]$.

(iii) *Assume that* (3.8) *and* (3.9) *hold. Then*

$$\lim_{n \to \infty} V_n^n(0) = \sup_{\nu \in \mathcal{U}(\mathbb{F})} \mathbb{E}[Y^\nu(T)].$$

PROOF. In view of Theorem 3.1, statement (iii) is a direct consequence of (i) and (ii). To see that (i) holds, we fix $\nu \in \mathcal{U}(\mathbb{F})$ and let $F^n$ denote the cumulative distribution of $T_n^n$. Let $\eta > 0$ be an arbitrary parameter. From the continuity condition (3.8), it follows that $|\mathbb{E}[Y^\nu(t)] - \mathbb{E}[Y^\nu(T)]| \leq \eta$ for $|T - t| \leq \varepsilon$ for sufficiently small $\varepsilon > 0$. Then, using Fubini's lemma together with the fact that the process $Y^\nu$ is bounded from below by a uniformly integrable martingale, it follows that

$$\mathbb{E}[Y^\nu(T_n^n)] \geq -C \mathbb{P}[|T_n^n - T| > \varepsilon] + \int_{T-\varepsilon}^{T+\varepsilon} \mathbb{E}[Y^\nu(t)] \, dF^n(t)$$

$$\geq -C \mathbb{P}[|T_n^n - T| > \varepsilon] + (\mathbb{E}[Y^\nu(T)] - \eta) \mathbb{P}[|T_n^n - T| \leq \varepsilon],$$

for some real constant $C > 0$. Since $T_n^n \longrightarrow T$ in probability, we deduce (i) by sending $n$ to $\infty$ and then $\eta$ to zero. Statement (ii) is obtained by following the lines of the above arguments, using the continuity condition (3.9). $\square$

REMARK 3.5. (i) The continuity assumptions (3.8), (3.9) have to be checked for each particular case; see, for example, Sections 4 and 5.

(ii) If there is some optimal control $\hat{\nu} \in \mathcal{U}(\mathbb{F})$ for the problem $\sup_{\nu \in \mathcal{U}(\mathbb{F})} \mathbb{E}[Y^\nu(T)]$, then it suffices to check condition (3.8) for $\hat{\nu}$.

(iii) The above proof provides an upper bound for the rate of convergence of $V_n^n$. Given the uniform modulus of continuity at $T$:

$$\rho(\varepsilon) := \sup_{t \in [T-\varepsilon, T+\varepsilon]} \sup_{\nu \in \mathcal{U}(\mathbb{F})} |\mathbb{E}[Y^\nu(t)] - \mathbb{E}[Y^\nu(T)]|,$$

the above arguments indeed show that

$$\left| V_n^n(0) - \sup_{\nu \in \mathcal{U}(\mathbb{F})} \mathbb{E}[Y^\nu(T)] \right| \leq C \mathbb{P}[|T_n^n - T| > \varepsilon] + \mathbb{P}[|T_n^n - T| \leq \varepsilon] \rho(\varepsilon)$$

for some real constant $C > 0$. Depending on $\rho$ and $T_n^n$, we can then choose $\varepsilon$ according to $n$ so as to minimize the right-hand side quantity. In general, $\rho$ is not known precisely but it is often possible to provide an upper bound which can be plugged in to the above inequality.

3.3. *Exponential maturity randomization.* In this section we assume that $(\zeta_n^j)_{j \leq n}$ is a sequence of exponentially distributed random variables with parameter $\lambda_n > 0$, for each $n$. In this case, (3.4) can be written as

$$e^{-\lambda_n t} V_{k+1}^\nu(t) = \operatorname*{ess\,sup}_{\mu \in \mathcal{U}(\mathbb{F}; t, \nu)} \mathbb{E}\left[ \lambda_n \int_t^\infty \bar{V}_k^\mu(u) e^{-\lambda_n u} \, du \Big| \mathcal{F}_t \right], \qquad t \geq 0,$$

so that the problem (3.4) is reduced to a classical stochastic control problem; see Remark 3.3. In this context, it suffices to assume that the bifurcation property (H$\mathcal{U}$) holds at $\mathbb{F}$-stopping times to obtain the existence of a measurable aggregating supermartingale; see [5].

For sake of completeness, we provide an easy proof of this result in the case where assumptions (HY), (H$\mathcal{U}$) and (H$\mathcal{U}$)″ are combined with a lower semicontinuity condition on $\nu \mapsto \mathbb{E}[Y^\nu(t)]$. In this case, we can even find a cadlag aggregating supermartingale.

LEMMA 3.2. *Let assumptions* (HY)–(H$\mathcal{U}$) *hold, and suppose that* $\mathcal{U}(\mathbb{F})$ *satisfies assumption* (H$\mathcal{U}$)″ *of Remark* 3.1. *Assume further that* $Y^\nu$ *is a cadlag process for each* $\nu \in \mathcal{U}(\mathbb{F})$, *and*

(3.10)
$$\liminf_{k \to \infty} \mathbb{E}[Y^{\nu_k}(t)] \geq \mathbb{E}[Y^\nu(t)]$$

*whenever* $\mathbb{P}[\nu_k(t) \longrightarrow \nu(t), \forall t \geq 0] = 1$.

*Then, for each* $k \leq n$ *and* $\nu \in \mathcal{U}(\mathbb{F})$, *we can find a cadlag supermartingale which aggregates* $V_k^\nu$ *in the sense of assumption* (HV).

PROOF. Since $V_k^\mu = V_k^\nu$ on $[0, t)$ for each $\mu \in \mathcal{U}(\mathbb{F}; t, \nu)$, we introduce the process

$$M_{k+1}^\nu(t) := e^{-\lambda_n t} V_{k+1}^\nu(t) + \lambda_n \int_0^t \bar{V}_k^\nu(u) e^{-\lambda_n u} \, du = \operatorname*{ess\,sup}_{\mu \in \mathcal{U}(\mathbb{F}; t, \nu)} J_{k+1}^\mu(t),$$

where

$$J_{k+1}^\mu(t) := \mathbb{E}\left[\lambda_n \int_0^\infty \bar{V}_k^\mu(u) e^{-\lambda_n u} \, du \,\Big|\, \mathcal{F}_t\right].$$

We first show that the process $M_{k+1}^\nu$ is a supermartingale for all $\nu \in \mathcal{U}(\mathbb{F})$ and $k \geq 0$. Indeed, under assumption (H$\mathcal{U}$), the family $\{J_{k+1}^\mu, \mu \in \mathcal{U}(\mathbb{F}; t, \nu)\}$ is directed upward. Then $M_{k+1}^\nu(t) = \lim_{n \to \infty} \uparrow J_{k+1}^{\mu_n}(t)$ for some sequence $(\mu_n)_{n \geq 1} \subset \mathcal{U}(\mathbb{F}; t, \nu)$, and it follows from the monotone convergence theorem that for all $s \leq t$,

$$\mathbb{E}[M_{k+1}^\nu(t) | \mathcal{F}_s] = \lim_{n \to \infty} \uparrow \mathbb{E}[J_{k+1}^{\mu_n}(t) | \mathcal{F}_s] = \lim_{n \to \infty} \uparrow J_{k+1}^{\mu_n}(s)$$

$$\leq \operatorname*{ess\,sup}_{\mu \in \mathcal{U}(\mathbb{F}; s, \nu)} J_{k+1}^\mu(s) = M_{k+1}^\nu(s).$$

We now turn to the proof of the statement of the lemma. We only show that $M_1^\nu$ admits a cadlag modification, and that, for each $t \geq 0$,

(3.11)
$$\liminf_{k \to \infty} \mathbb{E}[M_1^{\nu_k}(t)] \geq \mathbb{E}[M_1^\nu(t)]$$

whenever $\mathbb{P}[\nu_k(t) \longrightarrow \nu(t), \ t \geq 0] = 1.$

The required result will then be obtained by an induction argument.

We first prove that the mapping $t \mapsto \mathbb{E}[M_1^\nu(t)]$ is right-continuous. Since $M_1^\nu$ is a supermartingale, this ensures that it admits a cadlag modification; see, for example, [4]. First observe that, by the same argument as in Lemma 3.1, it follows from Assumption (H$\mathcal{U}$) that

(3.12) $$\mathbb{E}[M_1^\nu(t)] = \sup_{\mu \in \mathcal{U}(\mathbb{F}; t, \nu)} \mathbb{E}\left[\lambda_n \int_0^\infty Y^\mu(u) e^{-\lambda_n u} \, du\right].$$

This implies that $\mathbb{E}[M_1^\nu(t)]$ is nonincreasing in $t$. Hence, it suffices to show that

(3.13) $$\lim_{s \searrow t} \mathbb{E}[M_1^\nu(s)] \geq \mathbb{E}[M_1^\nu(t)].$$

To see this, fix $\varepsilon > 0$ and let $\mu_\varepsilon \in \mathcal{U}(\mathbb{F}; t, \nu)$ be such that

(3.14) $$\sup_{\mu \in \mathcal{U}(\mathbb{F}; t, \nu)} \mathbb{E}\left[\lambda_n \int_0^\infty Y^\mu(u) e^{-\lambda_n u} \, du\right] \leq \mathbb{E}\left[\lambda_n \int_0^\infty Y^{\mu_\varepsilon}(u) e^{-\lambda_n u} \, du\right] + \varepsilon.$$

Let $(t_k)_{k\geq 1}$ be a sequence converging toward $t$, and such that $t_k > t$, and define, for each $k \geq 1$,

$$\mu_\varepsilon^k := \nu \mathbb{1}_{[0,t_k)} + \mu_\varepsilon \mathbb{1}_{[t_k,\infty)}.$$

By assumption $(H\mathcal{U})''$, $\mu_\varepsilon^k \in \mathcal{U}(\mathbb{F}; t_k, \nu)$, so that by (3.12)

$$\mathbb{E}[M_1^\nu(t_k)] \geq \mathbb{E}\left[\lambda_n \int_0^\infty Y^{\mu_\varepsilon^k}(u) e^{-\lambda_n u}\, du\right].$$

Since $\mu_\varepsilon^k \longrightarrow \mu_\varepsilon$ $\mathbb{P}$-a.s., it follows from (3.10), (3.12), (3.14), Fatou's lemma, Remark 3.2 and Fubini's lemma that

$$\lim_{k\to\infty} \mathbb{E}[M_1^\nu(t_k)] \geq \liminf_{k\to\infty} \int_0^\infty \lambda_n \mathbb{E}[Y^{\mu_\varepsilon^k}(u)] e^{-\lambda_n u}\, du$$
$$\geq \int_0^\infty \lambda_n \mathbb{E}[Y^{\mu_\varepsilon}(u)] e^{-\lambda_n u}\, du$$
$$\geq \mathbb{E}[M_1^\nu(t)] - \varepsilon.$$

Sending $\varepsilon$ to 0 then shows (3.13).

Property (3.11) is easily deduced from (3.10) and (3.12) by using Fatou's and Fubini's lemmas as above. $\square$

## 4. Application 1: optimal stopping.

4.1. *The general case.* We now show that the optimal stopping problem presented in Section 2 can be embedded in the framework studied in the previous section. Let $Z$ be an $\mathbb{F}$-adapted process. We assume that $Z$ is cadlag and bounded by a uniformly integrable martingale. The main object of this section is the optimal stopping problem:

$$\sup_{\tau \in \mathcal{T}_0(\mathbb{F})} \mathbb{E}[Z(\tau \wedge T)].$$

In order to embed this problem in the general framework of the previous section, we follow [5] and set

$$\nu_\tau(t) := \mathbb{1}_{\tau < t} \qquad \text{for each } \tau \in \mathcal{T}_0(\mathbb{F}).$$

This defines a one-to-one correspondence between the set of stopping times $\mathcal{T}_0(\mathbb{F})$ and the family

$$\mathcal{U}(\mathbb{F}) := \{\nu_\tau : \tau \in \mathcal{T}_0(\mathbb{F})\}.$$

We shall denote by $\tau_\nu$ the stopping time associated to $\nu \in \mathcal{U}(\mathbb{F})$. Observing that

$$Z(\tau \wedge t) = Y^{\nu_\tau}(t) := \int_0^t Z(u)\, d\nu_\tau(u) + Z(t)\mathbb{1}_{\nu_\tau(t)=0},$$

we see that the optimal stopping problem can be rewritten as

$$\sup_{\tau \in \mathcal{T}_0(\mathbb{F})} \mathbb{E}[Z(\tau \wedge T)] = \sup_{\nu \in \mathcal{U}(\mathbb{F})} \mathbb{E}[Y^\nu(T)]. \tag{4.1}$$

REMARK 4.1. The set $\mathcal{U}(\mathbb{F})$ satisfies assumption (H$\mathcal{U}$) of Section 3. Also, for $\nu \in \mathcal{U}(\mathbb{F})$, $t \geq 0$ and $\mu \in \mathcal{U}(\mathbb{F}; t, \nu)$, we have $\tau_\mu = \tau_\nu$ on $\{\tau_\nu < t\}$. On $\{\tau_\nu \geq t\}$, $\tau_\mu$ can be any stopping time $\tau \in \mathcal{T}_0(\mathbb{F})$. However, assumption (H$\mathcal{U}$)″ is clearly not satisfied.

Given a sequence $(\zeta_n^k)_{k \leq n}$, we let $V_k^{\nu,n}$ be the associated sequence of controlled processes as defined in Section 3. Then, (3.3) reads as

$$V_0^{\nu,n}(t) = Z(\tau_\nu \wedge t), \qquad t \geq 0, \tag{4.2}$$

and it follows from Remark 4.1 that

$$V_1^{\nu,n}(t) = \operatorname*{ess-sup}_{\mu \in \mathcal{U}(\mathbb{F}; t, \nu)} \mathbb{E}[V_0^{\mu,n}(t + \zeta_n^n) | \mathcal{F}_t]$$

$$= \operatorname*{ess-sup}_{\mu \in \mathcal{U}(\mathbb{F}; t, \nu)} \mathbb{E}[Z(\tau_\mu \wedge (t + \zeta_n^n)) | \mathcal{F}_t]$$

$$= Z(\tau_\nu) \mathbb{1}_{\tau_\nu < t} + X_1^n(t) \mathbb{1}_{\tau_\nu \geq t},$$

where

$$X_1^n(t) := \operatorname*{ess-sup}_{\tau \in \mathcal{T}_t(\mathbb{F})} \mathbb{E}[Z(\tau \wedge (t + \zeta_n^n)) | \mathcal{F}_t], \qquad t \geq 0,$$

does not depend on $\tau_\nu$. We next compute

$$V_2^{\nu,n}(t) = \operatorname*{ess-sup}_{\mu \in \mathcal{U}(\mathbb{F}; t, \nu)} \mathbb{E}[\bar{V}_1^{\mu,n}(t + \zeta_n^{n-1}) | \mathcal{F}_t]$$

$$= \operatorname*{ess-sup}_{\mu \in \mathcal{U}(\mathbb{F}; t, \nu)} \mathbb{E}[Z(\tau_\nu) \mathbb{1}_{\tau_\nu < t} + \bar{X}_1^n(t + \zeta_n^{n-1}) \mathbb{1}_{\tau_\nu \geq t} | \mathcal{F}_t]$$

$$= Z(\tau_\nu) \mathbb{1}_{\tau_\nu < t} + X_2^n(t) \mathbb{1}_{\tau_\nu \geq t},$$

where

$$X_2^n(t) := \operatorname*{ess-sup}_{\tau \in \mathcal{T}_t(\mathbb{F})} \mathbb{E}[Z(\tau) \mathbb{1}_{\tau < t + \zeta_n^{n-1}} + \bar{X}_1^n(t + \zeta_n^{n-1}) \mathbb{1}_{\tau \geq t + \zeta_n^{n-1}} | \mathcal{F}_t], \qquad t \geq 0,$$

and $\bar{X}_1^n$ denotes a measurable aggregating process $X_1^n$ which we assume to exist. More generally, given $X_0^n := Z$ and

$$X_{k+1}^n(t) := \operatorname*{ess-sup}_{\tau \in \mathcal{T}_t(\mathbb{F})} \mathbb{E}[Z(\tau) \mathbb{1}_{\tau < t + \zeta_n^{n-k}} + \bar{X}_k^n(t + \zeta_n^{n-k}) \mathbb{1}_{\tau \geq t + \zeta_n^{n-k}} | \mathcal{F}_t],$$
(4.3)
$$t \geq 0$$

for $0 \leq k \leq n-1$, we observe the following relation between $V_k^{\nu,n}$ and $X_k^n$:

(4.4) $$V_k^{\nu,n}(t) = Z(\tau_\nu)\mathbb{1}_{\tau_\nu < t} + X_k^n(t)\mathbb{1}_{\tau_\nu \geq t}, \qquad t \geq 0.$$

In particular,

(4.5) $$V_n^n(0) = X_n^n(0),$$

and the sequence $(X_k^n(0))_{k \leq n}$ corresponds to Carr's algorithm as described in Section 2, for a suitable choice of $Z$.

We conclude this section with the following result which provides sufficient conditions for the convergence of the algorithm.

PROPOSITION 4.1. *Assume that $Z$ is cadlag and that assumption* (HY) *holds. Then,*

(4.6) $$\limsup_{\varepsilon \searrow 0} \sup_{\tau \in \mathcal{T}_0(\mathbb{F})} \mathbb{E}[|Z(\tau \wedge T) - Z(\tau \wedge (T+\varepsilon))|\mathbb{1}_{T<\tau}] = 0.$$

*In particular, if assumption* (HV) *holds and $T_n^n \to T$ in probability, then*

$$X_n^n(0) \longrightarrow \sup_{\tau \in \mathcal{T}_0(\mathbb{F})} \mathbb{E}[Z(\tau \wedge T)] \qquad \text{as } n \longrightarrow \infty.$$

PROOF. In view of (4.1)–(4.5), the second assertion is equivalent to

$$V_n^n(0) \longrightarrow \sup_{\nu \in \mathcal{U}(\mathbb{F})} \mathbb{E}[Y^\nu(T)] \qquad \text{as } n \longrightarrow \infty.$$

Observe that (4.6) implies (3.8), (3.9) of Corollary 3.1, so that the latter convergence result follows from Corollary 3.1(iii). It remains to show that (4.6) holds. For $\varepsilon > 0$, let $\tau^\varepsilon \in \mathcal{T}_0(\mathbb{F})$ be such that

$$\sup_{\tau \in \mathcal{T}_0(\mathbb{F})} \mathbb{E}[|Z(T) - Z(\tau \wedge (T+\varepsilon))|\mathbb{1}_{T<\tau}]$$
$$\leq \mathbb{E}[|Z(T) - Z(\tau^\varepsilon \wedge (T+\varepsilon))|\mathbb{1}_{T<\tau^\varepsilon}] + \varepsilon.$$

Since $Z$ is right-continuous,

$$\limsup_{\varepsilon \searrow 0} |Z(T) - Z(\tau^\varepsilon \wedge (T+\varepsilon))|\mathbb{1}_{T<\tau^\varepsilon} = 0, \qquad \mathbb{P}\text{-a.s.}$$

By the uniform integrability condition on $Z$, which is implied by assumption (HY), we deduce that

$$|Z(T) - Z(\tau^\varepsilon \wedge (T+\varepsilon))| \leq 2\sup_{t \geq 0}|Z(t)| \in L^1.$$

In view of the previous equality, the result follows from the dominated convergence theorem. □

4.2. *The case of exponentially distributed random variables.* In this section we discuss the case where, for each $n$, $(\zeta_n^j)_{j \leq n}$ is a sequence of exponentially distributed random variables with parameter $\lambda_n > 0$. Then, (4.3) can be written as

$$\underset{\tau \in \mathcal{T}_t(\mathbb{F})}{\text{ess-sup}} \mathbb{E}\left[Z(\tau)e^{-\lambda_n \tau} + \lambda_n \int_0^\tau \bar{X}_k^n(u)e^{-\lambda_n u}\, du \Big| \mathcal{F}_t\right]$$

(4.7)
$$= e^{-\lambda_n t} X_{k+1}^n(t) + \lambda_n \int_0^t \bar{X}_k^n(u)e^{-\lambda_n u}\, du.$$

In the case where $Z$ is cadlag and satisfies assumption (HY), we easily check that (HV) holds. In view of (4.4), this is implied by the next result.

LEMMA 4.1. *Assume that $Z$ is cadlag and that assumption* (HY) *holds. Then, for each $n \geq k \geq 1$, $X_k^n$ admits a cadlag aggregating supermartingale.*

PROOF. Assuming that $\bar{X}_k^n$ is of class (D), we deduce that the process

$$J_k^n(\cdot) := Z(\cdot)e^{-\lambda_n \cdot} + \lambda_n \int_0^\cdot \bar{X}_k^n(u)e^{-\lambda_n u}\, du$$

is of class (D) too. By Propositions 2.26 and 2.29 in [5], we deduce that the family

$$M_k^n(t) := \underset{\tau \in \mathcal{T}_t(\mathbb{F})}{\text{ess-sup}} \mathbb{E}[J(\tau)|\mathcal{F}_t]$$

can be aggregated by a supermartingale which is of class (D). The result then follows from (4.7) by induction. $\square$

In [2], the author considers the case where $Z(t) = e^{-rt}g(S(t))$, $t \geq 0$, for some function $g$, and a lognormal process $S$ defined by

$$S(t) = S(0)\exp\left[\left(r - \frac{\sigma^2}{2}\right)t + \sigma W(t)\right], \qquad t \geq 0,$$

for some real constants $r$, $\sigma$ and a standard Brownian motion $W$. It is shown that there is a sequence $(v_k^n)_{k \leq n}$ of bounded Lipschitz functions such that, for each $k \leq n$,

$$X_k^n = v_k^n(S).$$

Here, $X_k^n$ depends on time only through $S$. This is due to the time homogeneity of the dynamics of $S$.

For $g$ with polynomial growth and $\lambda_n = n$, it is clear that the conditions of Proposition 4.1 are satisfied for this simple model. This shows the consistency of Carr's algorithm.

## 5. Application 2: hedging in the uncertain volatility model.

5.1. *Problem formulation.* Let $W$ be a real-valued Brownian motion, on the probability space $(\Omega, \mathcal{F}, \mathbb{P})$, and let $\mathbb{F}$ be the $\mathbb{P}$-completion of the associated canonical filtration.

Given two constants $0 < \sigma_1 < \sigma_2$, we define $\mathcal{U}(\mathbb{F})$ as the collection of all $\mathbb{F}$-predictable processes $\nu$ with

$$\sigma_1 \leq \nu(\cdot) \leq \sigma_2, \qquad \mathbb{P}\text{-a.s.} \tag{5.1}$$

For each control process $\nu \in \mathcal{U}$, the controlled state process dynamics is defined by

$$dX^\nu(t) = X^\nu(t)\nu(t)\,dW(t), \qquad t \geq 0. \tag{5.2}$$

In this section we apply the maturity randomization technique to the stochastic control problem

$$v(0, x) := \sup_{\nu \in \mathcal{U}(\mathbb{F})} \mathbb{E}[h(X^\nu(T))|X^\nu(0) = x] \qquad \text{where } h : \mathbb{R}_+ \longrightarrow \mathbb{R} \tag{5.3}$$

is some bounded function. Further conditions will be placed later on $h$ in order to obtain an explicit maturity randomization algorithm.

The financial motivation of this problem is the following. The process $X^\nu$ represents the price of some given risky asset at each time $t$. $\nu$ is called the volatility process of $X^\nu$ and is only known to be bounded by two constants $\sigma_1$ and $\sigma_2$. The financial market also contains a nonrisky asset with price process normalized to unity. The random variable $h(X^\nu(T))$ is an example of European contingent claims. Then, $v(0, X^\nu(0))$ is the sharpest upper bound of all selling prices which are consistent with the no-arbitrage condition. We refer the readers to [11] and [8] for a deeper presentation of the theory of pricing contingent claims in general models. When $h$ is replaced by some convex (resp. concave) function, it was shown by El Karoui, Jeanblanc and Shreve [6] that the optimal control is $\nu^* \equiv \sigma_1$ (resp. $\nu^* \equiv \sigma_2$), and the associated hedging strategy is defined by the classical Black–Scholes strategy. The above simple model was studied by Avellaneda, Levy and Paras [1]. The connection with the hedging problem was analyzed by Cvitanić, Pham and Touzi [3] in the context of stochastic volatility models.

As usual, we introduce the dynamic version of the stochastic control problem (5.3):

$$v(t, x) := \sup_{\nu \in \mathcal{U}(\mathbb{F})} \mathbb{E}[h(X^\nu(T))|X^\nu(t) = x]. \tag{5.4}$$

Then, it follows from classical techniques that the function $v$ is the unique bounded $C^0([0, T] \times \mathbb{R}_+)$ viscosity solution of the nonlinear partial differential equation

$$-v_t - \tfrac{1}{2}x^2\sigma_2^2 v_{xx}^+ + \tfrac{1}{2}s^2\sigma_1^2 v_{xx}^- = 0 \qquad \text{on } [0, T) \times [0, \infty), v(T, \cdot) = h;$$

see, for example, [12]. Here subscripts denote partial derivatives. In the present context, we shall consider a function $h$ which is neither convex nor concave, so that no explicit solution for this PDE is known.

REMARK 5.1. Although more regularity should be expected for the value function $v$ because of the uniform parabolicity of the PDE, we do not enter this discussion since we only need the continuity property.

5.2. *Maturity randomization.* Let $(\xi^k)_{k \geq 0}$ be a sequence of independent random variables in $\mathcal{I}^{\mathbb{F}}$ with exponential distribution

$$\mathbb{P}[\xi^k \leq t] = 1 - e^{-t} \qquad \text{for all } k \geq 1.$$

Set

$$\zeta_n^k := \frac{T}{n} \xi^k \qquad \text{for every } k \leq n$$

so that

$$\sum_{k=1}^n \zeta_n^k \longrightarrow T, \qquad \mathbb{P}\text{-a.s.}$$

In the present context, the maturity randomization algorithm (3.3)–(3.4) translates to the sequence of stochastic control problems

$$v_n^0(x) = h(x)$$

and, for all $k \leq n-1$:

$$v_n^{k+1}(x) := \sup_{\nu \in \mathcal{U}(\mathbb{F})} \mathbb{E}[v_n^k(X^\nu(\zeta_n^{n-k}))|X^\nu(t) = x]$$

$$= \sup_{\nu \in \mathcal{U}(\mathbb{F})} \mathbb{E}\left[\int_0^\infty v_n^k(X^\nu(t)) \lambda_n e^{-\lambda_n t}\, dt \bigg| X^\nu(t) = x\right],$$

where $\lambda_n := n/T$. The corresponding Hamilton–Jacobi–Bellman equation is given by the ordinary differential equation (ODE)

(5.5) $\quad -\frac{1}{2} x^2 \sigma_2^2 [(v_n^{k+1})_{xx}]^+ + \frac{1}{2} x^2 \sigma_1^2 [(v_n^{k+1})_{xx}]^- + \lambda_n (v_n^{k+1} - v_n^k) = 0.$

An immediate induction argument shows that for each $1 \leq k \leq n$

(5.6) $\quad v_n^k$ is nonnegative, bounded, and satisfies $v_n^k(0) = 0$,

which provides the appropriate boundary condition for the above ODE.

We conclude this section by discussing the convergence of the maturity randomizing algorithm in this context, that is,

(5.7) $$\lim_{n \to \infty} v_n^n(X(0)) = v(0, X(0)).$$

Let $(V_n^{\nu,k})$ be defined as in Section 3:

$$V_0^{\nu,n} = h(X^\nu),$$

$$V_{k+1}^{\nu,n}(t) = \operatorname*{ess-sup}_{\mu \in \mathcal{U}(\mathbb{F};t,\nu)} \mathbb{E}[\bar{V}_k^{\mu,n}(t + \zeta_k^{n-k})|\mathcal{F}_t]$$

$$= \operatorname*{ess-sup}_{\mu \in \mathcal{U}(\mathbb{F};t,\nu)} \lambda_n \mathbb{E}\left[\int_t^\infty \bar{V}_k^{\mu,n}(u) e^{-\lambda_n(u-t)}\, du \Big| \mathcal{F}_t\right], \qquad t \geq 0, k \leq n-1,$$

so that, by the Markov feature of $X^\nu$,

$$V_k^{\nu,n} = v_n^k(X^\nu), \qquad 1 \leq k \leq n.$$

Clearly, assumption (HY) holds since $h$ is bounded. The above identity shows that assumption (HV) holds too.

We finally discuss conditions (3.8) and (3.9):

1. If $h$ is continuous, one deduces the a.s. continuity of $t \mapsto h(X^\nu(t))$ by using the bounds (5.1). Since $h$ is bounded, it follows that $t \mapsto \mathbb{E}[h(X^\nu(t))]$ is continuous too, that is, (3.8) holds.
2a. In the case where $h$ is Lipschitz continuous, (3.9) is deduced from the bounds of (5.1) and standard $L^2$ estimates on the diffusion process.
2b. In the case where $h$ is not Lipschitz continuous, we can use the fact that, as already noticed at the end of Section 5.1, the value function $v$ defined in (5.4) is continuous on $[0,T) \times (0,\infty)$. Since

$$v(\varepsilon,x) = \sup_{\nu \in \mathcal{U}(\mathbb{F})} \mathbb{E}[h(X^\nu(T-\varepsilon))|X^\nu(0)=x] \qquad \text{for } 0 < \varepsilon < T,$$

it follows from the homogeneity of the process $X^\nu$ that

$$t \mapsto \sup_{\nu \in \mathcal{U}(\mathbb{F})} \mathbb{E}[h(X^\nu(t))|X^\nu(0)=x]$$

is continuous.

5.3. *Explicit solution of the infinite horizon problems.* In this section we fix $n \geq 1$ and derive an explicit formula for the value function $v_n^{k+1}$ in terms of $v_n^k$ when the payoff function $h$ satisfies the following conditions:

(5.8) $\qquad\qquad\qquad\qquad h$ is continuous,

(5.9) $\quad h(x) = 1 - h(x^{-1}) = 0; \qquad 0 < x \leq x_0$ for some $x_0 \in (0,1)$,

and

(5.10) $\qquad\qquad h$ is convex on $[0,b_0]$, concave on $[b_0,\infty)$

$\qquad\qquad\qquad\qquad\qquad\qquad\qquad$ for some $x_0 < b_0 < x_0^{-1}$.

Notice that the above conditions imply that $h$ is nondecreasing on $\mathbb{R}_+$.

In order to derive an explicit expression of $v_n^{k+1}$ in terms of $v_n^k$, we shall exhibit a smooth solution $U^{k+1}$ of (5.5) which satisfies the properties (5.6). We then show that $U^{k+1} = v_n^{k+1}$ by a classical verification argument.

In view of the particular form of the function $h$, a bounded solution $U^{k+1}$ of (5.5) satisfying $U^{k+1}(0) = 0$ will be obtained under the additional guess that

$$(5.11) \qquad U_{xx}^{k+1}(x) \geq 0 \quad \text{if and only if } x \leq b_{k+1},$$

for some $b_{k+1} > 0$ to be determined. Then, the ODE (5.5) reduces to

$$(5.12) \qquad -\frac{1}{2}x^2\sigma_2^2 U_{xx}^{k+1} + \frac{n}{T}(U^{k+1} - U^k) = 0 \quad \text{for } x \leq b_{k+1},$$

$$(5.13) \qquad -\frac{1}{2}x^2\sigma_1^2 U_{xx}^{k+1} + \frac{n}{T}(U^{k+1} - U^k) = 0 \quad \text{for } x > b_{k+1}.$$

The solutions of (5.12) and (5.13) can be characterized by solving the associated homogeneous equations, and then applying the constants variation technique. Bounded solutions of (5.12) and (5.13) are then seen to be of the form

$$(5.14) \qquad U^{k+1}(x) = \begin{cases} A_1^{k+1}(x) x^{\gamma_1}, & \text{for } x > b_{k+1}, \\ A_2^{k+1}(x) x^{\gamma_2}, & \text{for } x \leq b_{k+1}, \end{cases}$$

where

$$\gamma_1 := \frac{1}{2}\left(1 - \sqrt{1 + \frac{8n}{T\sigma_1^2}}\right) \quad \text{and} \quad \gamma_2 := \frac{1}{2}\left(1 + \sqrt{1 + \frac{8n}{T\sigma_2^2}}\right).$$

We now plug (5.14) into (5.12)–(5.13). After some calculations, this leads to

$$A_i^{k+1}(x) = \gamma_i(1-\gamma_i)\int_{b_{k+1}}^{x} r^{-2\gamma_i} \int_{x_i}^{r} A_i^k(s) s^{2(\gamma_i-1)} \, ds \, dr$$

$$+ \alpha_i^{k+1} x^{1-2\gamma_i} + \beta_i^{k+1}, \qquad i = 1, 2,$$

where $x_i \geq 0$, $\alpha_i, \beta_i^{k+1}$, $i = 1, 2$, are constants to be fixed later on. By (5.14), this provides the candidate solution of (5.12)–(5.13):

$$U^{k+1}(x) = \begin{cases} \left(\dfrac{x}{b_{k+1}}\right)^{\gamma_1}\left[\beta_1^{k+1} + H_{b_{k+1}}^1[U^k]\left(\dfrac{x}{b_{k+1}}\right)\right] + \alpha_1^{k+1}\left(\dfrac{x}{b_{k+1}}\right)^{1-\gamma_1}, \\ \hfill x > b_{k+1}, \\ \left(\dfrac{x}{b_{k+1}}\right)^{\gamma_2}\left[\beta_2^{k+1} + H_{b_{k+1}}^2[U^k]\left(\dfrac{x}{b_{k+1}}\right)\right] + \alpha_2^{k+1}\left(\dfrac{x}{b_{k+1}}\right)^{1-\gamma_2}, \\ \hfill x \leq b_{k+1}, \end{cases}$$

where, for a function $\varphi : \mathbb{R}_+ \longrightarrow \mathbb{R}$, we denote

$$(5.15) \qquad H_b^i[\varphi](x) := \gamma_i(1-\gamma_i)\int_1^x r^{-2\gamma_i} \int_{x_i}^r \varphi(bs) s^{\gamma_i-2} \, ds \, dr.$$

In order to determine the constants $x_i, \alpha_i^{k+1}, \beta_i^{k+1}, i=1,2$, we now impose the restrictions of boundedness and nullity at zero:

$$(5.16) \quad \limsup_{x \nearrow \infty} U^{k+1}(x) = \limsup_{x \nearrow \infty} H^1_{b_{k+1}}[U^k](x) x^{\gamma_1} + \alpha_1^{k+1} x^{1-\gamma_1} < \infty,$$

$$(5.17) \quad \lim_{s \searrow 0} U^{k+1}(s) = \lim_{x \searrow 0} H^2_{b_{k+1}}[U^k] x^{\gamma_2} + \alpha_2^{k+1} x^{1-\gamma_2} = 0,$$

the continuity condition at the point $x = b_{k+1}$:

$$(5.18) \quad \beta_1^{k+1} + \alpha_1^{k+1} = \beta_2^{k+1} + \alpha_2^{k+1},$$

and the smooth-fit condition at the point $x = b_{k+1}$:

$$(5.19) \quad \begin{aligned} &\beta_1^{k+1}\gamma_1 + \{H^1_{b_{k+1}}[U^k]\}'(1) + \alpha_1^{k+1}(1-\gamma_1) \\ &= \beta_2^{k+1}\gamma_2 + \{H^2_{b_{k+1}}[U^k]\}'(1) + \alpha_2^{k+1}(1-\gamma_2). \end{aligned}$$

Since $1 - \gamma_1 > 0$ and $1 - \gamma_2 < 0$, it follows from the boundedness of $U^k$ that $H^i_{b_{k+1}}[U^k]$ is well defined with

$$(5.20) \quad x_1 = \infty \quad \text{and} \quad x_2 = 0.$$

We then conclude from (5.16)–(5.17) that

$$(5.21) \quad \alpha_1^{k+1} = \alpha_2^{k+1} = 0.$$

Using (5.18) and (5.19), it follows that

$$(5.22) \quad \beta_1^{k+1} = \beta_2^{k+1} = \beta[U^k](b_{k+1}) \quad \text{where } \beta[\varphi](b) := \int_0^\infty \varphi(br) f(r)\, dr$$

for any bounded function $\varphi : \mathbb{R}_+ \longrightarrow \mathbb{R}$, and

$$f(r) = \frac{1}{\gamma_2 - \gamma_1}[\gamma_2(\gamma_2-1) r^{\gamma_2-2} \mathbb{1}_{0 \le r \le 1} + \gamma_1(\gamma_1-1) r^{\gamma_1-2} \mathbb{1}_{r>1}].$$

For later use, we observe that

$$f > 0 \quad \text{on } (0, \infty) \quad \text{and} \quad \int_0^\infty f(r)\, dr = 1,$$

so that $f$ is a density function. In view of these results, we introduce the following notation. For a function $\varphi$ and some real constant $b > 0$, we set

$$(5.23) \quad T_b[\varphi](x) := \begin{cases} \left(\dfrac{x}{b}\right)^{\gamma_1} \left[\beta[\varphi](b) + H^1_b[\varphi]\left(\dfrac{x}{b}\right)\right], & \text{for } x > b, \\ \left(\dfrac{x}{b}\right)^{\gamma_2} \left[\beta[\varphi](b) + H^2_b[\varphi]\left(\dfrac{x}{b}\right)\right], & \text{for } x \le b, \end{cases}$$

so that our candidate solution can be written in the compact form

$$(5.24) \quad U^{k+1} = T_{b_{k+1}}[U^k] \quad \text{for some } b_{k+1} > 0.$$

REMARK 5.2. Let $(U^k)$ be a sequence defined as above with $U^0 = h$ satisfying (5.8), (5.9) and (5.10). As already observed, $U^0$ is nondecreasing and therefore nonnegative. As it is positive on some open set, one easily checks that $U^k(x) > 0$ for all $x > 0$ and $k \geq 1$ by using an inductive argument. Indeed, if $U^k$ is nonnegative, then $H^i_{b_{k+1}}[U^k] \geq 0$, $i = 1, 2$. If it is also positive on an open set, then $\beta[U^k](b_{k+1}) > 0$ whenever $b_{k+1} > 0$.

In order to fix the parameters $b_{k+1}$, we observe that if $U^{k+1}$ is convex on $[0, b_{k+1}]$ and concave on $[b_{k+1}, \infty)$, then it follows from (5.12)–(5.13) that $U^{k+1}(b_{k+1}) = U^k(b_{k+1})$. In view of (5.22), this provides the additional equation:

$$\beta[U^k](b_{k+1}) = U^k(b_{k+1}).$$

Our next results show that this condition defines uniquely the sequence of positive parameters $b_k$.

LEMMA 5.1. *Let $\varphi : \mathbb{R}_+ \longrightarrow [0,1]$ be a function satisfying*

$$(5.25) \quad \varphi(x) \sim_\infty 1 - a_1 x^{\gamma_1} (\ln x)^{\delta_1} \quad \text{and} \quad \varphi(x) \sim_0 a_2 x^{\gamma_2} (\ln x)^{\delta_2},$$

*for some positive constants $a_1, a_2$ and some integer $\delta_1, \delta_2$. Then there is a positive solution to the equation $\beta[\varphi](b) = \varphi(b)$, and*

$$T_b[\varphi](x) \sim_\infty 1 - a'_1 x^{\gamma_1} (\ln x)^{\delta'_1} \quad \text{and} \quad T_b[\varphi](x) \sim_0 a'_2 x^{\gamma_2} (\ln x)^{\delta'_2},$$

*for some positive constants $a'_1, a'_2$ and some integer $\delta'_1, \delta'_2$.*

PROOF. By the expression of the density $f$, it follows from a trivial change of variable that

$$(5.26) \quad \begin{aligned} \beta[\varphi](b) &= \frac{\gamma_2(\gamma_2 - 1)}{\gamma_2 - \gamma_1} b^{1-\gamma_2} \int_0^b r^{\gamma_2 - 2} \varphi(r) \, dr \\ &\quad + \frac{\gamma_1(\gamma_1 - 1)}{\gamma_2 - \gamma_1} b^{1-\gamma_1} \int_b^\infty r^{\gamma_1 - 2} \varphi(r) \, dr. \end{aligned}$$

Using the estimates of the lemma, we then compute that

$$\begin{aligned} \beta[\varphi](b) &\sim_0 \frac{\gamma_2(\gamma_2 - 1)}{\gamma_2 - \gamma_1} b^{1-\gamma_2} \int_0^b r^{\gamma_2 - 2} a_2 r^{\gamma_2} (\ln r)^{\delta_2} \, dr \\ &\quad + \frac{\gamma_1(\gamma_1 - 1)}{\gamma_2 - \gamma_1} b^{1-\gamma_1} \int_b^c r^{\gamma_1 - 2} a_2 r^{\gamma_2} (\ln r)^{\delta_2} \, dr + O(b^{1-\gamma_1}) \\ &\sim_0 a_2 b^{\gamma_2} (\ln b)^{\delta_2} \left[ \frac{\gamma_2(\gamma_2 - 1)}{(\gamma_2 - \gamma_1)(2\gamma_2 - 1)} \right. \end{aligned}$$

$$-\frac{\gamma_1(\gamma_1-1)}{(\gamma_2-\gamma_1)(\gamma_1+\gamma_2-1)}\Bigg] + O(b^{1-\gamma_1})$$

$$\sim_0 a_2 b^{\gamma_2}(\ln b)^{\delta_2}\left[1+\frac{\gamma_2(\gamma_2-1)}{(2\gamma_2-1)(1-\gamma_1-\gamma_2)}\right],$$

where the last equivalence follows from the fact that $\gamma_2 < 1 - \gamma_1$. From this, we conclude that

$$(5.27) \qquad \lim_{b\searrow 0}\frac{\beta[\varphi](b)}{\varphi(b)} = 1 + \frac{\gamma_2(\gamma_2-1)}{(2\gamma_2-1)(1-\gamma_1-\gamma_2)} > 1.$$

Next, since $f$ is a density, we have

$$1 - \beta[\varphi](b) = \frac{\gamma_2(\gamma_2-1)}{\gamma_2-\gamma_1} b^{1-\gamma_2} \int_0^b r^{\gamma_2-2}[1-\varphi(r)]\,dr$$
$$+ \frac{\gamma_1(\gamma_1-1)}{\gamma_2-\gamma_1} b^{1-\gamma_1} \int_b^\infty r^{\gamma_1-2}[1-\varphi(r)]\,dr.$$

By similar calculations, it follows from the estimate of the lemma that

$$(5.28) \qquad \lim_{b\nearrow\infty}\frac{1-\beta[\varphi](b)}{1-\varphi(b)} = \infty.$$

Now recall that $\varphi$ is continuous and bounded. Then $\beta[\varphi]$ is continuous, and the existence of a positive solution to the equation $\beta[\varphi](b) = \varphi(b)$ follows from (5.27) and (5.28).

The estimates on $T_b[\varphi]$ are deduced from (5.25) by similar arguments. $\square$

REMARK 5.3. The statement of Lemma 5.1 is valid for $\varphi = h$. Indeed, one can check that the above existence argument goes through under the condition (5.9) instead of (5.25).

LEMMA 5.2. *Let $\varphi:\mathbb{R}_+ \longrightarrow [0,1]$ be a nondecreasing function satisfying*

$$(5.29) \qquad \varphi(0) = 1 - \varphi(\infty) = 0$$

*such that*

(5.30) $\varphi$ *is convex on* $[0,b^*]$, *concave on* $[b^*,\infty)$ \qquad *for some $b^* > 0$*

*and either*:

(i) *there is some $\varepsilon > 0$ such that $\varphi(b) = 0$ for all $b \le \varepsilon$, or*
(ii) $\varphi$ *is strictly convex on a neighborhood of* 0.

*Then, there is at most one positive solution to the equation $\beta[\varphi](b) = \varphi(b)$.*

PROOF. Observe from (5.26) that the function $\beta[\varphi]$ is differentiable. From the convexity/concavity condition on $\varphi$, it follows that $\varphi$ is differentiable a.e. on $\mathbb{R}_+$, its subgradient $\partial_-\varphi$ is nonempty (resp. empty) in the domain of convexity (resp. concavity), and its supergradient $\partial_+\varphi$ is empty (resp. nonempty) in the domain of convexity (resp. concavity). Set $\partial\varphi := \partial_-\varphi \cup \partial_+\varphi$.

In order to prove the required result, it suffices to show that

$$\text{for all } b > 0 : \beta[\varphi](b) = \varphi(b) \implies \nabla\beta[\varphi](b) - p < 0 \tag{5.31}$$

for any $p \in \partial\varphi(b)$.

Recall that $\varphi(0) = 1 - \varphi(\infty) = 0$ by (5.29), and that $\varphi$ is nondecreasing, continuous on $[b^*, \infty)$. Since $f$ is density, it follows from (5.22) that $\beta[\varphi](b) > 0$ whenever $b > 0$, and therefore

$$0 = \varphi(0) < \varphi(b) < \varphi(\infty) = 1 \qquad \text{whenever } \beta[\varphi](b) = \varphi(b) \text{ with } b > 0. \tag{5.32}$$

With the help of (5.26), we next compute that

$$\nabla\beta[\varphi](b) = b^{-1}\gamma_1\varphi(b) + \gamma_1(\gamma_1 - 1)b^{-\gamma_1}\int_b^\infty \varphi(r)r^{\gamma_1 - 2}\,dr$$
$$+ b^{-1}(1 - \gamma_2)(\beta[\varphi] - \varphi)(b).$$

Integrating by parts the integral on the right-hand side, we see that

$$\nabla\beta[\varphi](b) = -\gamma_1 b^{-\gamma_1}\int_b^\infty \varphi'(r)r^{\gamma_1 - 1}\,dr + b^{-1}(1 - \gamma_2)(\beta[\varphi] - \varphi)(b),$$

so that

$$\nabla\beta[\varphi](b) = -\gamma_1 b^{-\gamma_1}\int_b^\infty \varphi'(r)r^{\gamma_1 - 1}\,dr \qquad \text{whenever } \beta[\varphi](b) = \varphi(b). \tag{5.33}$$

Similar computations provide the following alternative expression of the gradient:

$$\nabla\beta[\varphi](b) = \gamma_2 b^{-\gamma_2}\int_0^b \varphi'(r)r^{\gamma_2 - 1}\,dr \qquad \text{whenever } \beta[\varphi](b) = \varphi(b). \tag{5.34}$$

We now consider two cases:

1. Suppose that $b \geq b^*$ and choose an arbitrary $p \in \partial\varphi(b)$. The fact that $\varphi$ is concave nondecreasing on $[b, \infty)$ implies that $0 \leq \varphi'(r) \leq p$ for a.e. $r \geq b$. If $\varphi'(r) = p$ for a.e. $r \geq b$, we end up with a contradiction to (5.32). Hence, there is a subset of $[b, \infty)$ of positive measure on which $\varphi'(r) < p$ a.e. Together with (5.33) and the fact that $\gamma_1 < 0$, this implies that

$$\nabla\beta[\varphi](b) < -\gamma_1 b^{-\gamma_1} p \int_b^\infty r^{\gamma_1 - 1}\,dr = p \qquad \text{for any } p \in \partial\varphi(b).$$

Hence (5.31) holds in this case.

2. If $b \leq b^*$, we repeat the same argument as in the first case using the representation (5.34), and we show that (5.31) also holds in this case. □

We are now in a position to define our candidate solution of the nonlinear ODE (5.5).

PROPOSITION 5.1. *There exists a sequence of functions $(U^k)_{0 \leq k \leq n}$ defined by*

(5.35) $$U^0 = h \quad and \quad U^{k+1} = T_{b_{k+1}}[U^k],$$

*where the sequence $(b_k)_{k \geq 1}$ is uniquely defined by*

(5.36) $$\beta[U^k](b_{k+1}) = U^k(b_{k+1}),$$

*so that $U^{k+1}$ solves* (5.12)–(5.13). *Moreover, for all $k \geq 1$:*

(i) $U^k$ *is strictly convex (resp. strictly concave) on $(0, b_k)$ [resp. $(b_k, \infty)$],*
(ii) $(b_k - x)(U^k - U^{k-1})(x) > 0$ *for all $x \in (0, \infty) \setminus \{b_k\}$,*
(iii) $U^k$ *is a strictly increasing $C^2$ function with values in $[0, 1)$,*
(iv) $U^k(x) \sim_\infty 1 - a_1^k x^{\gamma_1} (\ln x)^{\delta_1^k}$ *and* $U^k(x) \sim_0 a_2^k x^{\gamma_2} (\ln x)^{\delta_2^k}$, *for some constants $a_1^k, a_2^k$ and some integer $\delta_1^k, \delta_2^k$.*

PROOF. 1. The existence and uniqueness of the sequence $(U^k)$ associated to the sequence $(b_k)$ follow from Lemmas 5.1 and 5.2, together with Remark 5.3, by a direct induction argument. The fact that $U^{k+1}$ solves (5.12)–(5.13) has been shown in the discussion preceding Lemma 5.1.

2. Item (iv) is then obtained by induction as a by-product of Lemma 5.1 and Remark 5.3. In view of (5.12)–(5.13), item (i) is a direct consequence of (ii). Also, $U^k$ is $C^2$ by construction, and the remaining part of item (iii) follows from (i), (iv) and an induction.

3. It only remains to prove item (ii). Clearly, it is sufficent to show that, for any $x > 0$,

(5.37) $$\{U^k\}'(x) - p < 0 \quad \text{for all } p \in \partial U^{k-1}(x)$$
$$\text{whenever } U^k(x) = U^{k-1}(x),$$

where we use the notation of Lemma 5.2. Indeed, this implies that $U^{k-1}$ and $U^k$ intersect at a unique point, which is already known to be $b_k$, and the required inequality follows. The reason for introducing the notation $\partial U^{k-1}(x)$ comes from the fact that, for $k = 1$, $U^0 = h$ may be nonsmooth although $h'$ is defined a.e. by (5.10). Let $x > 0$ be such that $U^k(x) = U^{k-1}(x)$ and set $i := 2$ if $x \leq b_k$ and $i := 1$ otherwise. From the expression of $U^k$ in terms of

$U^{k-1}$, we directly compute that

$$\{U^k\}'(x) = \frac{\gamma_i}{x} U^{k-1}(x) + x^{-\gamma_i} \gamma_i (1 - \gamma_i) \int_{x_i}^x U^{k-1}(r) r^{\gamma_i - 2} \, dr$$

$$= \frac{\gamma_i}{x} [U^k(x) - U^{k-1}(x)] + x^{-\gamma_i} \gamma_i \int_{x_i}^x \{U^{k-1}\}'(r) r^{\gamma_i - 1} \, dr$$

$$= x^{-\gamma_i} \gamma_i \int_{x_i}^x \{U^{k-1}\}'(r) r^{\gamma_i - 1} \, dr$$

by first integrating by parts and then using the assumption $U^k(x) = U^{k-1}(x)$.

3a. We first assume that $x \leq b_k$, so that the above identity reads

(5.38) $$\{U^k\}'(x) = x^{-\gamma_2} \gamma_2 \int_0^x \{U^{k-1}\}'(r) r^{\gamma_2 - 1} \, dr.$$

Fix $p \in \partial U^{k-1}(x)$. If $x \leq b_{k-1}$, we deduce from the convexity of $U^{k-1}$ on $[0, b_{k-1}]$ that $(U^{k-1})'(r) \leq p$ for a.e. $r \leq x$. Since $U^{k-1}(0) = 0$ and $x > 0$ implies $U^k(x) > 0$ by Remark 5.2, it follows from the nondecreasing feature of $U^{k-1}$, see (iii) and the remark just after (5.10), that $\{U^{k-1}\}'(r) < p$ a.e. on a subset of $[0, x]$ of positive measure. As $\gamma_2 > 0$, we deduce from (5.38) that

$$\{U^k\}'(x) < p$$

which is the required result.

If $x \in (b_{k-1}, b_k]$, then (5.38) can be written as

$$\{U^k\}'(x) = x^{-\gamma_2} \gamma_2 \int_0^{b_k} \{U^{k-1}\}'(r) r^{\gamma_2 - 1} \, dr - x^{-\gamma_2} \gamma_2 \int_x^{b_k} \{U^{k-1}\}'(r) r^{\gamma_2 - 1} \, dr.$$

By (5.36) and the identity (5.34) derived in the proof of Lemma 5.2, we obtain

$$\{U^k\}'(x) = \left(\frac{x}{b_k}\right)^{-\gamma_2} \nabla \beta [U^{k-1}](b_k) - x^{-\gamma_2} \gamma_2 \int_x^{b_k} \{U^{k-1}\}'(r) r^{\gamma_2 - 1} \, dr.$$

Since $x \geq b_{k-1}$, we deduce from the concavity of $U^{k-1}$ on $[b_{k-1}, \infty) \ni b_k$ that for $\hat{p} \in \partial U^{k-1}(b_k)$

$$\{U^k\}'(x) \leq \left(\frac{x}{b_k}\right)^{-\gamma_2} \nabla \beta [U^{k-1}](b_k) - x^{-\gamma_2} \gamma_2 \hat{p} \int_x^{b_k} r^{\gamma_2 - 1} \, dr$$

$$= \left(\frac{x}{b_k}\right)^{-\gamma_2} \nabla \beta [U^{k-1}](b_k) - \hat{p} \left[\left(\frac{x}{b_k}\right)^{-\gamma_2} - 1\right]$$

$$\leq \left(\frac{x}{b_k}\right)^{-\gamma_2} [\nabla \beta [U^{k-1}](b_k) - \hat{p}] + p.$$

Recalling the assertion (5.31) which was derived in the proof of Lemma 5.2, we deduce that $\{U^k\}'(x) - p < 0$ which concludes the proof in the case $x \leq b_k$.

3b. The case $x > b_k$ is treated similarly. Equation (5.38) is replaced by

$$\{U^k\}'(x) = x^{-\gamma_1}\gamma_1 \int_\infty^x \{U^k\}'(r) r^{\gamma_1 - 1}\, dr$$

and we use (5.33) instead of (5.34). $\square$

Our final result shows that the sequence $(U^k)_{k \leq n}$ constructed in the above proposition corresponds to $(v_n^k)_{k \leq n}$.

PROPOSITION 5.2. *Let $(U^k)_{k \leq n}$ be the sequence of functions defined in Proposition 5.1. Then, for each $k \geq 1$, $U^k = v_n^k$.*

PROOF. Since $U^0 = h$, it suffices to show that for all $x > 0$ and $k \geq 0$

$$U^{k+1}(x) = \sup_{\nu \in \mathcal{U}(\mathbb{F})} \mathbb{E}\left[\int_0^\infty U^k(X^\nu(t)) \frac{n}{T} e^{-(n/T)t}\, dt \Big| X^\nu(0) = x\right].$$

Let $k \geq 0$ be fixed. We first deduce from Proposition 5.1 that $U^{k+1}$ is a classical solution of

$$-\frac{1}{2}x^2\sigma_2^2[U_{xx}^{k+1}]^+ + \frac{1}{2}x^2\sigma_1^2[U_{xx}^{k+1}]^- + \frac{n}{T}(U^{k+1} - U^k) = 0 \quad \text{on } [0, \infty).$$

Since $\sigma_1 < \sigma_2$, the above ODE can be written as

$$(5.39) \qquad \sup_{\sigma_1 \leq \nu \leq \sigma_2} \frac{1}{2}x^2\nu^2 U_{xx}^{k+1} + \frac{n}{T}U^k = \frac{n}{T}U^{k+1} \quad \text{on } [0, \infty).$$

Recalling from Proposition 5.1 that $U^{k+1}$ is $C^2$, we then deduce from Itô's lemma that, for all $x \geq 0$, $\nu \in \mathcal{U}(\mathbb{F})$ and all stopping time $\tau$,

$$U^{k+1}(x) \geq \mathbb{E}\left[e^{-(n/T)\tau}U^{k+1}(X^\nu(\tau)) + \int_0^\tau U^k(X^\nu(t))\frac{n}{T}e^{-(n/T)t}\, dt \Big| X^\nu(0) = x\right].$$

Since $U^k$ and $U^{k+1}$ are bounded, it follows from the dominated convergence theorem that

$$(5.40) \qquad U^{k+1}(x) \geq \mathbb{E}\left[\int_0^\infty U^k(X^\nu(t))\frac{n}{T}e^{-(n/T)t}\, dt \Big| X^\nu(0) = x\right]$$

$$\text{for all } \nu \in \mathcal{U}(\mathbb{F}).$$

On the other hand, we deduce from (5.39) and Itô's lemma that for $\hat{\nu} \in \mathcal{U}(\mathbb{F})$ defined by

$$\hat{\nu}_t = \sigma_1 \mathbb{1}_{U_{xx}^{k+1}(X^{\hat{\nu}}(t)) < 0} + \sigma_2 \mathbb{1}_{U_{xx}^{k+1}(X^{\hat{\nu}}(t)) \geq 0}, \qquad t \geq 0,$$

we have, for all stopping time $\tau$,

$$U^{k+1}(x) = \mathbb{E}\left[e^{-(n/T)\tau}U^{k+1}(X_\tau^{\hat{\nu}}) + \int_0^\tau U^k(X_t^{\hat{\nu}})\frac{n}{T}e^{-(n/T)t}\,dt\Big|X^{\hat{\nu}}(0) = x\right].$$

Since $U^{k+1}$ and $U^k$ are bounded, we obtain by sending $\tau \to \infty$ that

$$U^{k+1}(x) = \mathbb{E}\left[\int_0^\infty U^k(X^{\hat{\nu}}(t))\frac{n}{T}e^{-(n/T)t}\,dt\Big|X^{\hat{\nu}}(0) = x\right],$$

which, combined with (5.40), concludes the proof. □

REMARK 5.4. Condition (5.9) can be clearly relaxed by only assuming that the payoff function satisfies the estimates (5.25) at infinity and the origin. We refrained from starting with such conditions because the parameters $\gamma_1$ and $\gamma_2$ arise along our analysis. Similarly, all our analysis goes through under the condition that $h$ is bounded, not necessarily lying in the interval $[0,1]$.

REMARK 5.5. Throughout this example, we assumed that the payoff function $h$ is continuous, excluding some important example in finance. The only place where this assumption was used is the proof Lemma 5.1 and to derive the continuity properties (3.8)–(3.9) of Corollary 3.1. Notice that some cases where $h$ is not continuous can be handled. Consider, for instance, the *digital option* example:

$$h(x) := \mathbb{1}_{[1,\infty)}(x) \qquad \text{for all } x \geq 0.$$

We directly compute that

$$U^1(x) := T_1[h](x) = [1 - (1-\beta_1)x^{\gamma_1}]\mathbb{1}_{[1,\infty)}(x) + \beta_1 x^{\gamma_2}\mathbb{1}_{[0,1)}(x),$$

where

$$\beta_1 = \frac{\gamma_1}{\gamma_1 - \gamma_2}.$$

When $h$ is continuous and satisfies the requirements of Lemmas 5.1 and 5.2, the constant $b_1$ is the unique solution of the equation $U^1(b_1) = h(b_1)$. In the above case of the digital option, notice that $h(1) = 1$ and $U^1(b_1) = \beta_1 \neq 1$. In particular, $U^1$ is not a $C^2$ function in this case.

Clearly, the above function $U^1$ is a bounded smooth solution of boths ODEs (5.12) and (5.13), and satisfies property (i) of Proposition 5.1. Although $U^1$ is not $C^2$ at the point $b_1 = 1$, the proof of Proposition 5.2 is still valid under the above properties, since Itô's lemma holds for the function $U^1$. Hence $U^1 = v_n^1$.

Observe that $U^1$ satisfies (5.25) of Lemma 5.1, and therefore Propositions 5.1 and 5.2 can be applied to the sequence $(U^k)$ started from $k = 2$.

By the same reasoning as in 2b in the discussion at the end of Section 5.2, the mapping

$$t \mapsto \sup_{\nu \in \mathcal{U}(\mathbb{F})} \mathbb{E}[h(X^\nu(t))|X^\nu(0) = x]$$

is continuous. For $\nu \in \mathcal{U}(\mathbb{F})$, we have $\nu \geq \sigma_1 > 0$ $\mathbb{P}$-a.s., so that $X^\nu$ is uniformly elliptic. This implies that

$$t \mapsto \mathbb{E}[h(X^\nu(t))|X^\nu(0) = x] = \mathbb{P}[X^\nu(t) \geq 1]$$

is also continuous. Hence Conditions (3.8) and (3.9) of Corollary 3.1 hold for this case.

5.4. *A numerical example.* In this section we use the maturity randomization algorithm to approximate the value function $v$ defined in (5.3). We consider the same model as in Section 5.1 with

(5.41) $$\sigma_1 = 0 \quad \text{and} \quad h(x) = \mathbb{1}_{x \geq K}$$

for some real parameter $K > 0$. The reasons for considering this particular case are:

1. The value function $v$ can be computed explicitly, up to a simple numerical integration. This will allow us to test our numerical results.
2. Although $\sigma_1 = 0$, the reasoning of Section 5.3 is easily adapted to this context.

PROPOSITION 5.3. *In the context of* (5.41), *the value function $v$ is given by*

$$v(0, x) = w(0, x)$$
$$:= \mathbb{1}_{x < K}\left[\int_{-\infty}^{m(x)} e^{-2m(x)(m(x)-r)/(\sigma_2^2 T)} f_T(r)\, dr + F_T(m(x))\right] + \mathbb{1}_{x \geq K}$$

*where* $m(x) := \ln(K/x)$ *and*

$$f_T(r) := \frac{1}{\sigma_2\sqrt{2\pi T}} e^{(-1/(2\sigma_2^2 T))(r+(1/2)\sigma_2^2 T)^2} \quad \text{and} \quad F_T(x) := \int_x^\infty f_T(r)\, dr.$$

*Furthermore, for every $x \geq 0$, the optimal control associated to $v(0, x)$ is given by*

$$\hat{\nu}_x(t) = \sigma_2 \mathbb{1}_{t \leq \tau_x}, \qquad t \in [0, T],$$

*where*

$$\tau_x := T \wedge \inf\{t \geq 0 : -\tfrac{1}{2}\sigma_2^2 t + \sigma_2 W_t \geq \ln(K/x)\}.$$

TABLE 1

| | $K = 100, x = 95, T = 0.5$ | | | | | | $K = 100, x = 95, T = 1$ | | | | |
|---|---|---|---|---|---|---|---|---|---|---|---|
| $\sigma_2 \backslash n$ | 10 | 200 | 500 | 1000 | Exact | $\sigma_2 \backslash n$ | 10 | 200 | 500 | 1000 | Exact |
| 0.2 | 0.6884 | 0.6978 | 0.6981 | 0.6982 | 0.6982 | 0.2 | 0.7693 | 0.7763 | 0.7765 | 0.7766 | 0.7767 |
| 0.4 | 0.8279 | 0.8330 | 0.8332 | 0.8333 | 0.8333 | 0.4 | 0.8697 | 0.8734 | 0.8735 | 0.8735 | 0.8736 |
| 0.6 | 0.8754 | 0.8789 | 0.8790 | 0.8790 | 0.8791 | 0.6 | 0.9030 | 0.9055 | 0.9056 | 0.9056 | 0.9056 |

PROOF. Clearly $w$ is continuous on $[0, T] \times [0, \infty)$ and $C^{1,2}$ on $[0, T] \times [0, K]$. Then, standard arguments show that it satisfies

$$(5.42) \qquad -v_t - \tfrac{1}{2}\sigma_2^2 x^2 v_{xx} = 0 \qquad \text{on } [0,T] \times [0,K],$$

and satisfies the boundary conditions

$$(5.43) \qquad v(T, \cdot) = \mathbb{1}_{\cdot \geq K} \quad \text{and} \quad v(\cdot, K) = 1.$$

For $\nu \in \mathcal{U}(\mathbb{F})$, let $X_{t,x}^{\nu}$ be the solution of (5.2) with initial condition $X_{t,x}^{\nu}(t) = x$ at time $t$. Recalling the law of the maximum of a drifted Brownian conditionally to its terminal value (see, e.g., [7]), we obtain that

$$(5.44) \quad w(t,x) = \mathbb{E}[h(X_{t,x}^{\hat{\nu}_{t,x}}(\tau_{t,x}^{\tilde{\nu}})) | X_{t,x}^{\hat{\nu}_{t,x}}(t) = x] = \mathbb{P}\left[\max_{t \leq s \leq T} X_{t,x}^{\hat{\nu}_{t,x}}(s) \geq K\right],$$

where

$$\hat{\nu}_{t,x}(s) = \sigma_2 \mathbb{1}_{s \leq \tau_{t,x}^{\tilde{\nu}}} \qquad \text{with } \tilde{\nu}(s) = \sigma_2, s \in [t, T]$$

and, for $\nu \in \mathcal{U}(\mathbb{F})$,

$$\tau_{t,x}^{\nu} := \inf\{t \leq s \leq T : X_{t,x}^{\nu}(s) \geq K\} \wedge T.$$

It follows that $w$ is nonincreasing in $t$. Since it solves (5.42), it is convex and solves

$$(5.45) \qquad \min_{0 \leq \nu \leq \sigma_2} -w_t - \tfrac{1}{2}\nu^2 x^2 w_{xx} = 0 \qquad \text{on } [0,T] \times [0,K].$$

Fix $\nu \in \mathcal{U}(\mathbb{F})$ and observe that, by Itô's lemma, (5.45), (5.43) and definition of $\tau_{t,x}^{\nu}$,

$$w(t,x) \geq \mathbb{E}[w(\tau_{t,x}^{\nu}, X_{t,x}^{\nu}(\tau_{t,x}^{\nu}))] = \mathbb{E}[\mathbb{1}_{\tau_{t,x}^{\nu} \leq T}] \geq \mathbb{E}[h(X_{t,x}^{\nu}(T))].$$

In view of (5.44), this implies that $w = v$ and that the optimal strategy is given by $\hat{\nu}_{t,x}$. $\square$

Table 2

| $\sigma_2 \backslash n$ | 10 | 200 | 500 | 1000 | Exact |
|---|---|---|---|---|---|
| | | | $K=100, x=50, T=1$ | | |
| 0.4 | $5.8058 \times 10^{-2}$ | $5.7949 \times 10^{-2}$ | $5.7951 \times 10^{-2}$ | $5.7952 \times 10^{-2}$ | $5.7954 \times 10^{-2}$ |
| | | | $K=100, x=80, T=1$ | | |
| 0.4 | $6.9973 \times 10^{-2}$ | $6.9430 \times 10^{-2}$ | $6.9419 \times 10^{-2}$ | $6.9415 \times 10^{-2}$ | $6.9411 \times 10^{-2}$ |

We next define the sequence of randomized control problems $(v_n^k)$ as in Section 5.2. The associated sequence of ODEs is given by

$$\min_{0 \leq \nu \leq \sigma_2} -\tfrac{1}{2}\nu^2 x^2 (v_n^{k+1})_{xx} + \lambda_n(v_n^{k+1} - v_n^k) = 0 \qquad \text{on } [0, K]$$

with $v_n^0 = h$ and $v_n^k(x) = 1$ for $x \geq K$. A straightforward adaptation of the arguments of Section 5.3 then shows that $(v_n^k)_{k \leq n}$ is explicitly given by the inductive scheme

$$v_n^{k+1} = \left(\frac{x}{K}\right)^{\gamma_2} \left[1 + H_K^2[v_n^k]\left(\frac{x}{K}\right)\right] \qquad \text{on } [0, K]$$

where $\gamma_2$ and $H^2$ are defined as in Section 5.3 for the corresponding value of $n$.

Condition (5.34) of Corollary 3.1 holds by Proposition 5.3. By Remark 3.5(ii), it suffices to check (5.33) for the optimal control associated to $v(0, x)$. Since this optimal control does not depend on the time horizon $T$, this amounts to checking (5.34). Since assumptions (HY), (HV) and (H$\mathcal{U}$) are satisfied, the above scheme is consistent.

In Tables 1 and 2, we report numerical estimates of $v$ obtained by using the approximating sequence $(v_n^n)$. The "exact" values of $v$ have been computed by numerical integration of the formula reported in Proposition 5.3.

In Table 1 above, we fix the parameters $K$, $x$, and we explore the performance of the maturity randomization algorithm for various values of $T$ and $\sigma_2$. Our experiments show an excellent performance of the algorithm. Notice that we already obtain sharp estimates for a small value of $n = 10$.

We next fix the parameter $\sigma_2$, and vary the values of the parameters $x$ and $T$. We observe again, in Table 2, the algorithm shows an excellent performance even for small values of $n$.

B. BOUCHARD
LPMA
UNIVERSITÉ PARIS VI
AND
CREST
PARIS
FRANCE
E-MAIL: bouchard@ccr.jussieu.fr

N. EL KAROI
CMAP
ECOLE POLYTECHNIQUE
PARIS
FRANCE
E-MAIL: elkaroui@cmapx.polytechnique.fr

N. TOUZI
CREST
AND
CEREMADE
UNIVERSITÉ PARIS IX
PARIS
FRANCE
E-MAIL: touzi@ensae.fr